\magnification\magstephalf
\documentstyle{amsppt}

\hsize 5.72 truein
\vsize 7.9 truein
\hoffset .39 truein
\voffset .26 truein
\mathsurround 1.67pt
\parindent 20pt
\normalbaselineskip 13.8truept
\normalbaselines
\binoppenalty 10000
\relpenalty 10000
\csname nologo\endcsname 


\font\bc=cmb10
\font\tenbsy=cmbsy10

\catcode`\@=11

\def\myitem#1.{\item"(#1)."\advance\leftskip10pt\ignorespaces}

\def\qedsymbol{{\mathsurround\z@$\square$}}
\redefine\qed{\relaxnext@\ifmmode\let\next\@qed\else
  {\unskip\nobreak\hfil\penalty50\hskip2em\null\nobreak\hfil
    \qedsymbol\parfillskip\z@\finalhyphendemerits0\par}\fi\next}
\def\@qed#1$${\belowdisplayskip\z@\belowdisplayshortskip\z@
  \postdisplaypenalty\@M\relax#1
  $$\par{\lineskip\z@\baselineskip\z@\vbox to\z@{\vss\noindent\qed}}}
\outer\redefine\beginsection#1#2\par{\par\penalty-250\bigskip\vskip\parskip
  \leftline{\tenbsy x\bf#1. #2}\nobreak\smallskip\noindent}
\outer\redefine\genbeginsect#1\par{\par\penalty-250\bigskip\vskip\parskip
  \leftline{\bf#1}\nobreak\smallskip\noindent}

\def\next{\let\@sptoken= }\def\next@{ }\expandafter\next\next@
\def\@futureletnext#1{\let\nextii@#1\futurelet\next\@flti}
\def\@flti{\ifx\next\@sptoken\let\next@\@fltii\else\let\next@\nextii@\fi\next@}
\expandafter\def\expandafter\@fltii\next@{\futurelet\next\@flti}

\let\zeroindent\z@
\let\savedef@\endproclaim\let\endproclaim\relax 
\define\chkproclaim@{\add@missing\endroster\add@missing\enddefinition
  \add@missing\endproclaim
  \envir@stack\endproclaim
  \edef\endit@{\leftskip\the\leftskip\rightskip\the\rightskip}}
\let\endproclaim\savedef@
\def\thing@{.\enspace\egroup\ignorespaces}
\def\thingi@(#1){ \rm(#1)\thing@}
\def\thingii@\cite#1{ \rm\@pcite{#1}\thing@}
\def\thingiii@{\ifx\next(\let\next\thingi@
  \else\ifx\next\cite\let\next\thingii@\else\let\next\thing@\fi\fi\next}
\def\thing#1#2#3{\chkproclaim@
  \ifvmode \medbreak \else \par\nobreak\smallskip \fi
  \noindent\advance\leftskip#1
  \hskip-#1#3\bgroup\bc#2\unskip\@futureletnext\thingiii@}
\let\savedef@\endproclaim\let\endproclaim\relax 
\def\endit{\endproclaim\endit@\let\endit@\undefined}
\let\endproclaim\savedef@
\def\defn#1{\thing\parindent{Definition #1}\rm}
\def\lemma#1{\thing\parindent{Lemma #1}\sl}
\def\prop#1{\thing\parindent{Proposition #1}\sl}
\def\thm#1{\thing\parindent{Theorem #1}\sl}
\def\cor#1{\thing\parindent{Corollary #1}\sl}
\def\conj#1{\thing\parindent{Conjecture #1}\sl}

\def\remk#1{\thing\zeroindent{Remark #1}\rm}

\def\narrowthing#1{\chkproclaim@\medbreak\narrower\noindent
  \it\def\next{#1}\def\next@{}\ifx\next\next@\ignorespaces
  \else\bgroup\bc#1\unskip\let\next\narrowthing@\fi\next}
\def\narrowthing@{\@futureletnext\thingiii@}

\def\@cite#1,#2\end@{{\rm([\bf#1\rm],#2)}}
\def\cite#1{\in@,{#1}\ifin@\def\next{\@cite#1\end@}\else
  \relaxnext@{\rm[\bf#1\rm]}\fi\next}
\def\@pcite#1{\in@,{#1}\ifin@\def\next{\@cite#1\end@}\else
  \relaxnext@{\rm([\bf#1\rm])}\fi\next}

\advance\minaw@ 1.2\ex@
\atdef@[#1]{\ampersand@\let\@hook0\let\@twohead0\brack@i#1,\z@,}
\def\brack@{\z@}
\let\@@hook\brack@
\let\@@twohead\brack@
\def\brack@i#1,{\def\next{#1}\ifx\next\brack@
  \let\next\brack@ii
  \else \expandafter\ifx\csname @@#1\endcsname\brack@
    \expandafter\let\csname @#1\endcsname1\let\next\brack@i
    \else \Err@{Unrecognized option in @[}%
  \fi\fi\next}
\def\brack@ii{\futurelet\next\brack@iii}
\def\brack@iii{\ifx\next>\let\next\brack@gtr
  \else\ifx\next<\let\next\brack@less
    \else\relaxnext@\Err@{Only < or > may be used here}
  \fi\fi\next}
\def\brack@gtr>#1>#2>{\setboxz@h{$\m@th\ssize\;{#1}\;\;$}%
 \setbox@ne\hbox{$\m@th\ssize\;{#2}\;\;$}\setbox\tw@\hbox{$\m@th#2$}%
 \ifCD@\global\bigaw@\minCDaw@\else\global\bigaw@\minaw@\fi
 \ifdim\wdz@>\bigaw@\global\bigaw@\wdz@\fi
 \ifdim\wd@ne>\bigaw@\global\bigaw@\wd@ne\fi
 \ifCD@\enskip\fi
 \mathrel{\mathop{\hbox to\bigaw@{$\ifx\@hook1\lhook\mathrel{\mkern-9mu}\fi
  \setboxz@h{$\displaystyle-\m@th$}\ht\z@\z@
  \displaystyle\m@th\copy\z@\mkern-6mu\cleaders
  \hbox{$\displaystyle\mkern-2mu\box\z@\mkern-2mu$}\hfill
  \mkern-6mu\mathord\ifx\@twohead1\twoheadrightarrow\else\rightarrow\fi$}}%
 \ifdim\wd\tw@>\z@\limits^{#1}_{#2}\else\limits^{#1}\fi}%
 \ifCD@\enskip\fi\ampersand@}
\def\brack@less<#1<#2<{\setboxz@h{$\m@th\ssize\;\;{#1}\;$}%
 \setbox@ne\hbox{$\m@th\ssize\;\;{#2}\;$}\setbox\tw@\hbox{$\m@th#2$}%
 \ifCD@\global\bigaw@\minCDaw@\else\global\bigaw@\minaw@\fi
 \ifdim\wdz@>\bigaw@\global\bigaw@\wdz@\fi
 \ifdim\wd@ne>\bigaw@\global\bigaw@\wd@ne\fi
 \ifCD@\enskip\fi
 \mathrel{\mathop{\hbox to\bigaw@{$%
  \setboxz@h{$\displaystyle-\m@th$}\ht\z@\z@
  \displaystyle\m@th\mathord\ifx\@twohead1\twoheadleftarrow\else\leftarrow\fi
  \mkern-6mu\cleaders
  \hbox{$\displaystyle\mkern-2mu\copy\z@\mkern-2mu$}\hfill
  \mkern-6mu\box\z@\ifx\@hook1\mkern-9mu\rhook\fi$}}%
 \ifdim\wd\tw@>\z@\limits^{#1}_{#2}\else\limits^{#1}\fi}%
 \ifCD@\enskip\fi\ampersand@}


\define\ie{{\it i.e\.}}
\define\today{\number\day\ \ifcase\month\or
  January\or February\or March\or April\or May\or June\or
  July\or August\or September\or October\or November\or December\fi
  \ \number\year}
\def\pr@m@s{\ifx'\next\let\nxt\pr@@@s \else\ifx^\next\let\nxt\pr@@@t
  \else\let\nxt\egroup\fi\fi \nxt}

\define\widebar#1{\mathchoice
  {\setbox0\hbox{\mathsurround\z@$\displaystyle{#1}$}\dimen@.1\wd\z@
    \ifdim\wd\z@<.4em\relax \dimen@ -.16em\advance\dimen@.5\wd\z@ \fi
    \ifdim\wd\z@>2.5em\relax \dimen@.25em\relax \fi
    \kern\dimen@ \overline{\kern-\dimen@ \box0\kern-\dimen@}\kern\dimen@}%
  {\setbox0\hbox{\mathsurround\z@$\textstyle{#1}$}\dimen@.1\wd\z@
    \ifdim\wd\z@<.4em\relax \dimen@ -.16em\advance\dimen@.5\wd\z@ \fi
    \ifdim\wd\z@>2.5em\relax \dimen@.25em\relax \fi
    \kern\dimen@ \overline{\kern-\dimen@ \box0\kern-\dimen@}\kern\dimen@}%
  {\setbox0\hbox{\mathsurround\z@$\scriptstyle{#1}$}\dimen@.1\wd\z@
    \ifdim\wd\z@<.28em\relax \dimen@ -.112em\advance\dimen@.5\wd\z@ \fi
    \ifdim\wd\z@>1.75em\relax \dimen@.175em\relax \fi
    \kern\dimen@ \overline{\kern-\dimen@ \box0\kern-\dimen@}\kern\dimen@}%
  {\setbox0\hbox{\mathsurround\z@$\scriptscriptstyle{#1}$}\dimen@.1\wd\z@
    \ifdim\wd\z@<.2em\relax \dimen@ -.08em\advance\dimen@.5\wd\z@ \fi
    \ifdim\wd\z@>1.25em\relax \dimen@.125em\relax \fi
    \kern\dimen@ \overline{\kern-\dimen@ \box0\kern-\dimen@}\kern\dimen@}%
  }

\catcode`\@\active
\input epsf

\font\tenscr=rsfs10 
\font\sevenscr=rsfs7 
\font\fivescr=rsfs5 
\skewchar\tenscr='177 \skewchar\sevenscr='177 \skewchar\fivescr='177
\newfam\scrfam \textfont\scrfam=\tenscr \scriptfont\scrfam=\sevenscr
\scriptscriptfont\scrfam=\fivescr
\define\scr#1{{\fam\scrfam#1}}
\let\Cal\scr

\catcode`\@=11
\def\orderatorname#1{\mathord{\newmcodes@\kern\z@\fam\z@#1}}
\catcode`\@\active

\let\0\relax 
\define\ConjA{A}
\define\ConjB{B}
\define\exc{{\text{exc}}}
\define\Gm{\Bbb G_{\text{m}}}
\define\Id{\operatorname{Id}}
\define\Image{\operatorname{Im}}
\define\ord{\operatorname{ord}}
\define\Pic{\operatorname{Pic}}
\define\pr{\orderatorname{pr}}

\define\Supp{\operatorname{Supp}}
\define\restrictedto#1{\big|_{#1}}

\topmatter
\title On the abc conjecture and diophantine approximation by rational points%
  \endtitle
\rightheadtext{On abc and approximation by rational points}
\author Paul Vojta\endauthor
\affil University of California, Berkeley\endaffil
\address Department of Mathematics, University of California,
  970 Evans Hall\quad\#3840, Berkeley, CA \ 94720-3840\endaddress
\date 9 December 1999\enddate
\thanks Supported by NSF grants DMS95-32018 and DMS99-70393,
the Institute for Advanced Study, and IHES.\endthanks

\abstract
We show that an earlier conjecture of the author, on diophantine approximation
of rational points on varieties, implies the ``abc conjecture'' of Masser
and Oesterl\'e.  In fact, a weak form of the former conjecture is sufficient,
involving an extra hypothesis that the variety and divisor admit a faithful
group action of a certain type.  Analogues of this weaker conjecture are
proved in the split function field case of characteristic zero, and in the
case of holomorphic curves (Nevanlinna theory).

The proof of the latter involves a geometric generalization of the
classical lemma on the logarithmic derivative, due to McQuillan.
This lemma may be of independent interest.
\endabstract
\endtopmatter

\document

This paper discusses some conjectures that, if true, would imply the
following conjecture, known as the Masser-Oesterl\'e ``abc conjecture.''

\conj{\00.1} (Masser-Oesterl\'e)  Let $\epsilon>0$.  Then there is a constant
$C$, depending only on $\epsilon$, such that for all triples $a,b,c\in\Bbb Z$
with $a+b+c=0$ and $(a,b,c)=1$, the following inequality holds:
$$\max\{|a|,|b|,|c|\} \le C\prod_{p\mid abc}p^{1+\epsilon}\;.\tag\00.1.1$$
\endit

It has been known for some time that this conjecture would follow from other
well-known conjectures; for example, see \cite{Vo~1, \S5.ABC}.  In particular,
it would be a consequence of the following conjecture:

\conj{\00.2} \cite{Vo~1, Conjecture 5.2.6}  Let $k$ be a number field, let $S$
be a finite set of places of $k$ containing all the archimedean places of $k$,
let $X$ be a smooth projective curve over $k$, let $D$ be an effective divisor
on $X$ without multiple points, let $\Cal K$ denote the canonical line sheaf
on $X$, let $\Cal A$ be an ample line sheaf on $X$, let $\epsilon>0$,
and let $r$ be a positive integer.  Then the inequality
$$m_{S,k}(D,P) + h_{\Cal K,k}(P)
  \le \frac{\log|D_{k(P)/\Bbb Q}|}{[k(P):k]}
    + \epsilon\,h_{\Cal A,k}(P) + O(1)$$
holds for all $P\in X(\widebar{\Bbb Q})$ with $P\notin\Supp D$
and $[k(P):k]\le r$.
\endit

Here $h_{\Cal K,k}$ and $h_{\Cal A,k}$ denote logarithmic heights normalized
relative to $k$, and $m_{S,k}$ is the proximity function for $D$.
See Section \01 for details.

Although it is stated here only for curves, this conjecture is still out of
reach at the present time, due to the problem of dealing with the discriminant
term.

The purpose of the present paper is to show how the abc conjecture would follow
from the following, possibly weaker, conjecture.

\conj{\ConjA}  \cite{Vo~1, Conjecture 3.4.3}  Let $k$ be a number field, let $S$
be a finite set of places of $k$ containing all the archimedean places of $k$,
let $X$ be a smooth complete variety over $k$, let $D$ be a normal crossings
divisor on $X$ (assumed effective and without multiple components),
let $\Cal K$ denote the canonical line sheaf on $X$, let $\Cal A$ be a
big line sheaf on $X$, and let $\epsilon>0$.  Then there exists a proper
Zariski-closed subset $Z$ of $X$, depending only on $X$, $D$, $\Cal A$, and
$\epsilon$, such that the inequality
$$m_{S,k}(D,P) + h_{\Cal K,k}(P) \le \epsilon\,h_{\Cal A,k}(P) + O(1)$$
holds for all $P\in \bigl(X\setminus(Z\cup\Supp D)\bigr)(k)$.
\endit

This conjecture has the obvious advantage of dealing only with rational points,
at the expense of allowing $X$ to have arbitrary dimension.  In fact,
to get arbitrarily small values of $\epsilon$, Conjecture \ConjA{} would
need to be known for certain pairs $(X,D)$ of arbitrarily large dimension.

On the other hand, the pairs $(X,D)$ are special in the sense that they admit
a faithful group action of $\Gm^{\dim X-1}$.  Therefore, it would be sufficient
to prove the following weakening of Conjecture \ConjA:

\conj{\ConjB}  Conjecture \ConjA{} holds under the additional assumption
that there is a semiabelian variety $G$, of dimension $\dim X-1$,
acting faithfully on $X$ in such a way that the action preserves $D$.
\endit

If all points in $X$ had finite stabilizers under the group action,
then a quotient $X/G$ might exist, and would be a curve.
In that case Conjecture \ConjB{} may possibly follow from the fact that
Conjecture \ConjA{} is known for curves.  Unfortunately, there always exist
points in $X$ with infinite stabilizers, so this approach does not work.
However, even without forming a quotient, it is possible to prove analogues of
Conjecture \ConjB{} in the split function field case of characteristic zero
and in the case of holomorphic curves $f\:\Bbb C\to X$ (Nevanlinna theory).

The proof of the latter involves a result of McQuillan.
An immediate corollary of this result, Corollary \05.2, gives a
Nevanlinna-like inequality involving pull-backs of holomorphic differential
forms on the domain space.  This corollary corresponds to the trivial fact
for algebraic maps, that if $f\:C\to X$ is an algebraic map from a complex
projective curve $C$ to a projective variety $X$ and $\omega$
is a meromorphic differential form on $X$ such that $f^{*}\omega\ne0$, then
the number of zeroes minus the number of poles of $f^{*}\omega$ must equal
$2g(C)-2$.  This corollary may be useful for translating theorems about
curves on varieties into corresponding theorems on holomorphic curves.

The proof of McQuillan's result mentioned above involves
a geometric generalization of the classical lemma on the logarithmic derivative,
Theorem \0A.2, which may also be of interest.  An earlier geometric
version of the logarithmic derivative lemma was proved by Noguchi
\cite{N, Lemma~2.3}.  His version was stated in terms of a (nonzero) global
section of $\Omega^1_X(\log D)$; it follows from the present version.  
Others have generalized Noguchi's work, but McQuillan was the first to
remove the dependence on global sections.

This paper is organized as follows.  Section \01 sets notation and recalls
some of the basic definitions.  Section \02 gives a characterization of the
exceptional set $Z$ in Conjecture \ConjB{}.
Section \03, which forms the heart of the paper, contains the proof that
Conjecture \ConjB{} implies the abc conjecture.

Sections \04 and \05 prove analogues of Conjecture \ConjB{} in the split
function field case and in the Nevanlinna theory case, respectively.
Section \06 introduces a hierarchy of variations of Conjecture \ConjA{}
and similar diophantine statements.  For example, it has been known for
decades that it is often productive to prove a statement in the function field
case before trying to prove it for number fields.  The hierarchy contains this,
as well as the corresponding observation about Nevanlinna theory.
In light of this hierarchy, Sections \04 and \05 represent the progression
of Conjecture \ConjB{} through the steps in the hierarchy.

Section \07 of the paper gives another way of obtaining a weak form
of the abc conjecture (for $\epsilon>26$) from Conjecture \ConjA{}.
This particular variant works with one particular three-fold, with $D=0$.
It also gives an example of how Conjecture \ConjB{} on a certain rational
projective surface would imply something abc-like.

Finally, the paper concludes with an appendix giving a proof of
Proposition \05.1, since a proof has not appeared in print.
This proof follows the ideas of McQuillan \cite{McQ~2}, where the $D=0$
case is proved, and of P.-M. Wong \cite{W}.  See also \cite{McQ~3}.

The author thanks William Cherry for help with questions in Nevanlinna theory.

\beginsection{\01}{Notation and definitions}

In this paper, a {\bc variety} is an integral scheme, separated and of finite
type over a field.  The conventions regarding $\Bbb P(\Cal E)$ and $\Cal O(1)$
on $\Bbb P(\Cal E)$, where $\Cal E$ is a vector sheaf, are as in EGA.
In other words, a point on $\Bbb P(\Cal E)$ corresponds to a hyperplane
in the corresponding fiber of $\Cal E$.

Unless otherwise specified, $k$ will denote a global field of
characteristic zero, $M_k$ its set of places, and $S$ a finite subset of $M_k$
containing all archimedean places.  Each place $v\in M_k$ has an associated
almost-absolute value $\|\cdot\|_v$, normalized as follows.

If $k$ is a number field with ring of integers $R$, then
$$\|x\|_v = \cases |\sigma(x)|^{[k_v:\Bbb R]} &\; \text{if $v$ is archimedean,
    corresponding to $\sigma,\bar\sigma\:k\hookrightarrow\Bbb C$;} \\
  (R:\frak p)^{-\ord_{\frak p}(x)} &\; \text{\vtop{\hbox{if $v$ is
    non-archimedean, corresponding to}
    \hbox{a prime ideal $(0)\ne\frak p\subseteq R$.}}}
  \endcases$$
If $v$ is a complex place, then $\|\cdot\|_v$ fails to satisfy the triangle
inequality, hence the term {\it almost\/}-absolute value.

If $k$ is a function field, then we assume without loss of generality that its
field of constants is algebraically closed, and define
$\|x\|_v=\exp(-\ord_v(x))$ for $x\in k^{\times}$.  (Of course, the abc
conjecture itself is already known in the function field case, but some
parts of this paper require the notation in the function field case.)

With this choice of normalization, the product formula
$\prod_{v\in M_k}\|x\|_v=1$, $x\in k^{\times}$, holds without multiplicities.

If $P\in\Bbb P^n(k)$ is a rational point with homogeneous coordinates
$[x_0:\dots:x_n]$, then we define the {\bc height}
$$h_k(P) = \sum_{v\in M_k}\log\max\{\|x_0\|_v,\dots,\|x_n\|_v\}\;.$$
(We shall always assume that homogeneous coordinates lie in the field
of definition of the point.)  If $E$ is a finite extension of $k$
(with compatible global field structure), then $h_E(P)=[E:k]h_k(P)$,
so we may define $h_k(P)$ for $P\in\Bbb P^n(\bar k)$ by $h_k(P)=h_E(P)/[E:k]$
for any field $E$ containing $k(P)$.

The well-known functoriality and additivity properties of Weil heights
then allow us to define a height $h_{\Cal L,k}\:X(\bar k)\to\Bbb R$
for any complete variety $X$ over $k$ and any line sheaf $\Cal L$ on $X$.
These heights are only defined up to $O(1)$.

For places $v\in M_k$ and Cartier divisors $D$ on a complete variety $X$, let
$$\lambda_{D,v}\:(X\setminus\Supp D)(\bar k_v)\to\Bbb R$$
be Weil functions, normalized so that if $D=(f)$ locally then $\lambda_{D,v}$
behaves like $-\log\|f\|_v$ near $\Supp D$.  Then the height function may
be decomposed as a sum
$$h_{\Cal O(D),k}(X) = \sum_{v\in M_k}\lambda_{D,v}(P) + O(1)\tag\01.1$$
for all $P\in X(k)\setminus\Supp D$.  See \cite{L, Chap.~10} for details on
Weil functions.

For a finite subset $S$ of $M_k$, Weil functions allow us to break the height
$h_{\Cal O(D),k}$ up into two parts, the {\bc proximity function}
$$m_{S,k}(D,P) := \sum_{v\in S}\lambda_{D,v}(P)$$
and the {\bc counting function}
$$N_{S,k}(D,P) := \sum_{v\notin S}\lambda_{D,v}(P)\;,$$
so that
$$h_{\Cal O(D),k}(P) = m_{S,k}(D,P) + N_{S,k}(D,P) + O(1)$$
for all $P\in X(k)\setminus\Supp D$, by (\01.1).  As with the height,
the proximity and counting functions satisfy
$$m_{T,E}(D,P)=[E:k]m_{S,k}(D,P) \qquad\text{and}\qquad
  N_{T,E}(D,P)=[E:k]N_{S,k}(D,P)\;,$$
where $T$ is the set of places in $M_E$ lying over places in $S$.  This
allows us to define $m_{S,k}$ and $N_{S,k}$ for $P\in X(\bar k)\setminus\Supp D$
so that (\01.2) still holds.  We also note that the proximity and counting
functions are additive and functorial in $D$.

\beginsection{\02}{The exceptional set in Conjecture \ConjB}

This section characterizes the exceptional set $Z$ in Conjecture \ConjB{}
in terms of the group action.

We start with some facts about groups acting on varieties.
Let $G$ be a group variety over a field $k$ of characteristic zero,
and let $X$ be a smooth variety over $k$ on which $G$ acts.
We have a morphism
$$\alpha\:G\times X\to X\times X$$
whose first component is the map $G\times X\to X$ defining the group action
and whose second component is the projection to $X$.  Viewing both as schemes
over $X$ via the second projection, the morphism $\alpha$ becomes a morphism
over $X$.  This induces a morphism of relative tangent bundles:
$$\alpha_{*}\:T_{G\times X/X}\to\alpha^{*}T_{X\times X/X}\;.$$
Pulling back the above map by the map $(0,\Id_X)\:X\to G\times X$
and using the isomorphisms $T_G\cong\Cal O_G^{\dim G}$ and
$T_{X\times X/X}\cong\pr_1^{*}T_X$ gives a map
$$\Cal O_X^{\dim G}\to(0,\Id_X)^{*}\alpha^{*}\pr_1^{*}T_X
  \cong T_X\;.$$
Taking $\wedge^{\dim G}$ of this map then gives a map
$$\Cal O_X \to \wedge^{\dim G}T_X\;,$$
which is equivalent to giving a section
$$\omega\in\Gamma(X,\wedge^{\dim G}T_X)\;.$$

\lemma{\02.1}  The above section $\omega$ is nonzero if and only if the kernel
of the group action has dimension zero.
\endit

\demo{Proof}  Suppose first that the kernel $H$ has nonzero dimension.
Since $H$ is normal, $G/H$ exists, and $\alpha$ factors through
$\pi\:G\times X\to(G/H)\times X$.  Therefore $\alpha_{*}$ factors through
$\pi_{*}\:T_{G\times X/X}\to\pi^{*}T_{(G/H)\times X/X}$.
Since $\wedge^{\dim G}T_{(G/H)}=0$, the section $\omega$ is zero.

Conversely, suppose that $\dim H=0$.  Let $\eta$ be the generic point of $X$,
and let $L=K(\eta)$.  The map $G_{\widebar L}\to X_{\widebar L}$ defined
by $g\mapsto g\eta$ has finite fibers, hence is \'etale.  This implies that
$\alpha_{*}$ is an isomorphism at $(0,\eta)$, so $\omega$ is nonzero
at $\eta$.\qed
\enddemo

\defn{\02.2}  Let $G$ be a group variety acting on a variety $X$.  An orbit
$Gx$ (for $x\in X$) is {\bc degenerate} if $\dim Gx<\dim G$.
\endit

\lemma{\02.3}  Let $G$ be a commutative group variety over a field $k$ of
characteristic zero, and let $X$ be a smooth variety over $k$ on which $G$
acts faithfully.  Then there exists an open dense $G$\snug-invariant subset
$U\subseteq X$, a variety $V$ over $k$, and a morphism $\pi\:U\to V$
exhibiting $V$ as a quotient $U/G$.
\endit

\demo{Proof}  Let $U_0$ be the open set on which the above-defined
section $\omega$ is nonzero.  Since $G$ is commutative, $\alpha$ is
$G$\snug-equivariant (if we let $G$ act on $G\times X$ by acting on the
second factor, and on $X\times X$ by acting on both factors); hence $\omega$
is $G$\snug-invariant, so $U_0$ is also $G$\snug-invariant.
Since $G$ acts faithfully, Lemma \02.1 implies that $\omega\ne0$,
so $U_0\ne\emptyset$.

Moreover, Lemma \02.1 applied to each orbit in $U_0$ implies that the stabilizer
of every point in $U_0$ is finite.  Therefore \cite{SGA~3, V~8.1} applies,
giving an open dense subset $U\subseteq U_0$ and a quotient morphism
$\pi\:U\to U/G$.\qed
\enddemo

\defn{\02.4}  Let $X$ be a scheme, let $P\in X$ be a regular point,
and let $D$ be a divisor on $X$.  Then $D$ has {\bc normal crossings}
at $P$ if, locally in the \'etale topology, $D$ is effective and can be
written as a principal divisor $D=(x_1\dotsm x_r)$, where $x_1,\dots,x_r$
are elements of the maximal ideal $\frak m_{P'}$ in the local ring
at the point $P'$ in the \'etale neighborhood,
and the images of $x_1,\dots,x_r$ in $\frak m_{P'}/\frak m_{P'}^2$
are linearly independent over the residue field at $P'$.
We say that $D$ is a {\bc normal crossings divisor} on a regular scheme $X$
if it has normal crossings at all $P\in X$.
\endit

Note that, under this definition, a normal crossings divisor must be
effective and reduced.  The notion of normal crossings does not make sense
at a singular point.

\remk{\02.5}  Definition \02.4 can be restated as follows.  A divisor $D$
has normal crossings at $P$ if and only if there exists an \'etale
neighborhood $\pi\:X'\to X$ of $P$ and functions $x_1,\dots,x_r\in\Cal O(X')$,
such that (i) $\pi^{*}D=(x_1\dotsm x_r)$, and (ii) for each subset
$I\subseteq\{1,\dots,r\}$, the subscheme cut out by the ideal generated
by $\{x_i:i\in I\}$ is regular of codimension $\#I$.
\endit

\defn{\02.6}  Let $D$ be a divisor on a scheme $X$.  Then we say that
the pair $(X,D)$ is {\bc regular} at a point $P\in X$ if $X$ is regular
at $P$ and $D$ has normal crossings there.  We say that the pair $(X,D)$
is {\bc regular}, or that $(X,D)$ is a {\bc regular pair}, if $(X,D)$ is
regular at all $P\in X$.
\endit

\defn{\02.7}  Let $(X,D)$ be a regular pair, and let $f\:X\to Y$ be a dominant
morphism.  We say that $(X,D)$ has {\bc good reduction} at a point $Q\in Y$
if the pair $(X_Q,D\restrictedto{X_Q})$ is regular, where $X_Q$ is the fiber
of $f$ over $Q$.
\endit

\prop{\02.8}  Let $f\:X\to Y$ be a dominant morphism of integral, separated
schemes of finite type over $\Bbb Z$, and let $D$ be a divisor on $X$.
If the pair $(X,D)$ is regular, then it has good reduction over an open
dense subset of $Y$.
\endit

\demo{Proof}  First note that $(X,D)$ has good reduction over the generic
point of $Y$.  Indeed, the condition of Definition \02.4 remains true when
restricting to the generic fiber.  Thus it will suffice to show that the set
$$\bigl\{P\in X:
  \text{$(X_{f(P)},D\restrictedto{X_{f(P)}})$ is regular at $P$}\bigr\}$$
is an open subset of $X$.  But this follows from the corresponding fact for
smooth morphisms, via Remark \02.5.\qed
\enddemo

\lemma{\02.9}  Let $i\:G\hookrightarrow\widebar G$ be an equivariant
completion of a semiabelian variety over an algebraically closed field
of characteristic zero, such that the pair $(\widebar G,\widebar G\setminus G)$
is regular.  Then the logarithmic canonical line sheaf
$\Cal K_{\widebar G}(\widebar G\setminus G)$ is trivial.
\endit

\demo{Proof}  This was proved in the first few sentences of the proof of
\cite{Vo~3, Lemma~5.6}, but under the additional assumption that the equivariant
completion is of the form indicated in \cite{Vo~3, Lemma~2.2}.  It remains to
be shown that all equivariant completions of semiabelian varieties are of
this form (\ie, obtained in a natural way from an equivariant completion
of the toric part).

Let $\rho\:G\to A$ be the maximal abelian quotient of $G$, and
let $k$ be the ground field.  Since $\widebar G$ is a nonsingular variety,
$\rho$ extends to a morphism $\bar\rho\:\widebar G\to A$ \cite{Mi, Thm.~3.1}.
Since $k$ is algebraically closed, the kernel of $\rho$ is a split torus,
so every point of $A$ has a Zariski-open neighborhood $U$ such that
$\rho^{-1}(U)$ is a product $U\times\Gm^\mu$.
In particular, $\rho$ has a rational section over $U$.  By equivariance,
$\bar\rho^{-1}(U)$ must therefore be of the form $\bar\rho^{-1}(0)\times U$,
hence $\widebar G$ is of the desired form.\qed
\enddemo

\prop{\02.10}  Let $k$, $S$, $X$, $D$, $\Cal A$, and $\epsilon$ be as in
Conjecture \ConjA.  Assume:
\roster
\myitem i.  a semiabelian variety $G$ of dimension $\dim X-1$ acts faithfully
on $X$, preserving $D$;
\myitem ii.  the support of $D$ contains all degenerate orbits;
\myitem iii.  the map $\pi\:U\to V$ of Lemma \02.3 extends to a morphism
$\bar\pi\:X\to\widebar V$ for some projective completion $\widebar V$ of $V$;
and
\myitem iv.  Conjecture \ConjA{} holds for the above data.
\endroster
Then Conjecture \ConjA{} holds for the above data with $Z=\bar\pi^{-1}(\Sigma)$,
where $\Sigma$ is the set of all points on $\widebar V$ over which
the pair $(X,D)$ has bad reduction.
\endit

\demo{Proof}  Let $Z$ be an exceptional set for which Conjecture \ConjA{} holds.

There exists a finite subset $\Gamma\subseteq G(k)$ such that
$$Z' := \bigcap_{g\in\Gamma}\bigl(T_g^{-1}(Z)\cup T_{-g}^{-1}(Z)\bigr)$$
is preserved under the action of $G$, where $T_g\:X\to X$ denotes
translation by $g$.  This can be rewritten
$$Z' = \bigcup_\sigma \bigcap_{g\in\Gamma} T_{\sigma(g)g}^{-1}(Z)\;,$$
where $\sigma$ ranges over all functions $\sigma\:\Gamma\to\{\pm1\}$.

Since the height and proximity functions are functorial, condition (iv)
implies that for any fixed $g\in G$,
$$m_{S,k}(T_g^{*}D, P) + h_{\Cal K,k}(P)
  \le \epsilon\,h_{T_g^{*}\Cal A,k}(P) + O(1)\tag\02.10.1$$
for all $P\in(X\setminus T_g^{*}(Z))(k)$ (with a possibly different constant
$O(1)$, depending on $g$).  Clearly $T_g^{*}D=D$. Moreover,
$$T_g^{*}\Cal A\otimes T_{-g}^{*}\Cal A \cong \Cal A^{\otimes 2}\;.$$
Indeed, this holds because $g\mapsto T_g^{*}\Cal A\otimes\Cal A\spcheck$
defines a morphism $G\to\Pic^0(X)$, which must be a group homomorphism.  Thus,
$$h_{\Cal A,k}(P)
  \ge \min\{h_{T_g^{*}\Cal A,k}(P),h_{T_{-g}^{*}\Cal A,k}(P)\} + O(1)$$
for all $P\in X(\bar k)$, where the constant $O(1)$ depends on $g$ but not
on $P$.  Consequently, we have
$$\split h_{\Cal A,k}(P)
  &\ge \max_{g\in\Gamma}\min\{h_{T_g^{*}\Cal A,k}(P),h_{T_{-g}^{*}\Cal A,k}(P)\}
    + O(1) \\
  &\ge \min_{\sigma\in\{\pm1\}^\Gamma} \max_{g\in\Gamma}
    h_{T_{\sigma(g)g}^{*}\Cal A,k}(P) \;. \endsplit$$
Thus (\02.10.1), applied for all $g$ such that $\pm g\in\Gamma$, implies that
Conjecture \ConjA{} holds with exceptional set $Z'$ (after adjusting $O(1)$).
Therefore we may assume that $Z$ is preserved by $G$.

Let $Y$ be an irreducible component of $Z$ not contained in $\Supp D$.
Then $Y$ is preserved by the group action, so it contains
a non-degenerate orbit.  For dimension reasons, it follows that $Y$ is the
closure of a non-degenerate orbit.  Thus $Y$ is an equivariant completion
of a semiabelian variety, and by the assumptions on $D$, it follows that
the semiabelian variety of which $Y$ is the closure, is the complement
of $Y\cap\Supp D$.

In addition, $Y$ is contained in a fiber of $\bar\pi$.  Indeed, let $G$ act
trivially on $V$ and on $\widebar V$.  Then $\pi$ is $G$\snug-equivariant,
so the same is true for $\bar\pi$.  Thus fibers of $\bar\pi$ are
$G$\snug-invariant.

Suppose now that $Y$ does not lie on a fiber of bad reduction for $(X,D)$.
Then, for dimension reasons (since $\widebar V$ is a curve) $Y$ is the entire
fiber of $\bar\pi$, and the pair $(Y,D\restrictedto Y)$ is regular since
the fiber has good reduction.  As mentioned above, $Y$ is an equivariant
completion of a semiabelian variety, so Lemma \02.9 implies that
the line sheaf $\Cal K_Y\bigl(D\restrictedto Y\bigr)$ is trivial.
It then follows easily that Conjecture \ConjA{} holds on $Y$
with empty exceptional set.  But also, since $Y$ is an entire fiber,
we have that $\Cal O(Y)\restrictedto Y$ is trivial,
so $\Cal K_Y\cong\Cal K_X\restrictedto Y$, and therefore Conjecture \ConjA{}
on $Y$ implies that all points on $Y$ satisfy Conjecture \ConjA{} on $X$
(after suitably adjusting the $\Cal O(1)$ term).  Thus $Y$ can be removed
from $Z$.\qed
\enddemo

\remk{\02.11}  It is not always possible to take $Z=\emptyset$ in
Conjecture \ConjB, as the following example illustrates.
Let $X$ be the blow-up of $\Bbb P^1\times\Bbb P^1$ at the point $(0,1)$,
and let $D$ be the pull-back of the divisor
$[0]\times\Bbb P^1+[\infty]\times\Bbb P^1+\Bbb P^1\times[0]
  +\Bbb P^1\times[\infty]$.
Let $\Gm$ act on $\Bbb P^1\times\Bbb P^1$ by acting in the obvious way on
the first factor; this extends to an action of $\Gm$ on $X$ that preserves $D$.
In this situation, Conjecture \ConjA{} does not holds for $(X,D)$ unless
$Z$ contains the strict transform $Y$ of $\Bbb P^1\times\{1\}$.
This is true because the restriction $\Cal K_X\restrictedto Y$ is isomorphic
to $\Cal O(-1)$ instead of $\Cal O(-2)$ (via the isomorphism $Y\cong\Bbb P^1$).
\endit

\remk{\02.12}  Conditions (ii) and (iii) do not necessarily restrict the
applicability of Proposition \02.10.  Given a regular pair $(X,D)$ with a group
action as in (i), one can find another regular pair $(X',D')$ satisfying
(i)--(iii), admitting a proper birational morphism $\pi\:X'\to X$ such that
$\Supp D'\supseteq\Supp\pi^{*}D$ and such that the group action extends
to $X'$ and preserves $D'$.  In that situation, Conjecture \ConjA{}
for $(X',D')$ implies Conjecture \ConjA{} for $(X,D)$, so if Conjecture \ConjA{}
is known for $(X',D')$, then the exceptional set for $(X,D)$ will be contained
in the image of the exceptional set for $(X',D')$.
\endit

\beginsection{\03}{Conjecture \ConjB{} and the abc conjecture}

This section shows that Conjecture \ConjB, if proved, would imply the abc
conjecture.

In this section, all heights, proximity functions, and counting functions
are taken relative to $S=\{\infty\}$ over $\Bbb Q$, so the
subscripts will be omitted from the notation.  The abc conjecture can also
be formulated over number fields, with arbitrary (finite) set $S$.
The methods of this sections extend readily to this case.

\defn{\03.1}  For $\ell\in\Bbb Z$, $\ell\ne0$, let
$$n(\ell) = \sum_{p\mid\ell}\log p$$
(where the sum is over distinct primes dividing $\ell$).
\endit

Then the main inequality of the abc conjecture can be replaced by
$$n(abc) \ge (1-\epsilon)h([a:b:c]) - O_\epsilon(1)\tag\03.2$$
(with a different $\epsilon$).  We will use this inequality in place of
(\00.1.1).

The basic idea of this section is to construct, for a given $\epsilon$,
a regular pair $(\Gamma',D')$, and for each triple $(a,b,c)$ as in the
abc conjecture a point $P_{a,b,c}'\in\Gamma'(\Bbb Q)$ such that
Conjecture \ConjB{} applied to the points $P_{a,b,c}'$ implies the
abc conjecture with the desired value of $\epsilon$.

To begin, fix $0<\epsilon<1$, and let $n$ be an integer with $n>3/\epsilon$.

For triples $a,b,c\in\Bbb Z$ with $a+b+c=0$ and $(a,b,c)=1$, let $x_1,\dots,x_n$
be integers such that
$$x_1x_2^2x_3^3\dotsm x_n^n = a$$
and such that
$$\ord_p(x_i) = \cases \fracwithdelims[]{\ord_p(a)}{n}&\text{if $i=n$;} \\
  1&\text{if $i=\ord_p(a)-n\fracwithdelims[]{\ord_p(a)}{n}$; and} \\
  0&\text{otherwise.}\endcases$$
Likewise, choose $y_1,\dots,y_n$ and $z_1,\dots,z_n$ such that
$$y_1y_2^2y_3^3\dotsm y_n^n=b
  \qquad\text{and}\qquad z_1z_2^2z_3^3\dotsm z_n^n=c\;.$$
This defines a point $P=P_{a,b,c}$ in $\bigl(\Bbb P^2\bigr)^n$ with
multihomogeneous coordinates
$$([x_1:y_1:z_1],[x_2:y_2:z_2],\dots,[x_n:y_n:z_n])\;.$$
All such points lie on the variety $X_n\subseteq\bigl(\Bbb P^2)^n$
defined by the equation
$$\prod_{i=1}^n x_i^i + \prod_{i=1}^n y_i^i + \prod_{i=1}^n z_i^i = 0\;.$$

We will want to compare the height of $P_{a,b,c}\in X_n$ with the
height $h([a:b:c])$ occurring in $(\03.2)$.  This can be done
via functoriality of heights, as follows.  The map
$$\bigl([x_1:y_1:z_1],\dots,[x_n:y_n:z_n]\bigr)\mapsto
  \biggl(\prod x_i^i:\prod y_i^i:\prod z_i^i\biggr)$$
defines a rational map $X_n\dashrightarrow\Bbb P^2$ (actually to the line
$a+b+c=0$ in $\Bbb P^2$).  Let $\Gamma_n$ be the closure of the graph of
this rational map, and let $\phi\:\Gamma_n\to\Bbb P^2$ be the corresponding
morphism.  Then $\Gamma_n\subseteq\bigl(\Bbb P^2\bigr)^{n+1}$ is defined
by the equations:
$$\align \prod x_i^i + \prod y_i^i + \prod z_i^i &= 0 \tag\03.3a \\
  a + b + c &= 0 \tag\03.3b \\
  a\prod y_i^i &= b\prod x_i^i\;.  \tag\03.3c \endalign$$
This system appears to be asymmetrical in the variables, but in fact we have
$$\align 0
  &= a\Bigl(\prod x_i^i + \prod y_i^i + \prod z_i^i\Bigr)
    - (a + b + c)\prod x_i^i \\
  &= \Bigl(a\prod y_i^i - b\prod x_i^i\Bigr)
    + \Bigl(a\prod z_i^i - c\prod x_i^i\Bigr) \endalign$$
and similarly for $b\prod z_i^i - c\prod y_i^i$.  Thus $\Gamma_n$ is preserved
by the symmetric group $\Cal S_3$ acting by permuting the variables.

Given a triple $(a,b,c)$ as above, let $P^{*}=P^{*}_{a,b,c}$ be the
(unique) point on $\Gamma_n$ lying over $P_{a,b,c}\in X_n$.  Then we have
$$h([a:b:c]) = h_{\phi^{*}\Cal O(1)}(P^{*}) + O(1)\;.\tag\03.4$$

From now on, let $D$ be the divisor $x_1\dotsm x_ny_1\dotsm y_nz_1\dotsm z_n=0$
on $\Gamma_n$, and let $E$ be the divisor $[1:-1:0]+[0:1:-1]+[-1:0:1]$
on the line $a+b+c=0$ in $\Bbb P^2$.  The next step is to compare $n(abc)$
with $N(D,P^{*}_{a,b,c})$.

\lemma{\03.5}  For relatively prime integers $a,b,c$ with $a+b+c=0$,
$$N(D,P^{*}_{a,b,c}) \le n(abc) + \frac1n N(E,\phi(P^{*}_{a,b,c})) + O(1)\;.
  \tag\03.5.1$$
\endit

\demo{Proof}  To begin, note that
$$\lambda_{E,v}([a:b:c])
  := -\log\frac{\|abc\|_v}{\max\{\|a\|_v,\|b\|_v,\|c\|_v\}^3}$$
and
$$\split &\lambda_{D,v}([a:b:c],[x_1:y_1:z_1],\dots,[x_n:y_n:z_n]) \\
  &\qquad := -\log\frac{\|x_1\dotsm x_ny_1\dotsm y_nz_1\dotsm z_n\|_v}
    {\max\{\|x_1\|_v,\|y_1\|_v,\|z_1\|_v\}^3\dotsm
      \max\{\|x_n\|_v,\|y_n\|_v,\|z_n\|_v\}^3}\endsplit$$
are Weil functions for $E$ on $a+b+c=0$ in $\Bbb P^2$ and for $D$ on $X$,
respectively.  We assume here that $a,b,c$ are relatively prime and that
$x_i,y_i,z_i$ are relatively prime for all $i$.  Under that assumption,
the above expressions simplify to
$$\lambda_{E,v}([a:b:c]) = \ord_p(abc)\log p$$
and
$$\lambda_{D,v}([a:b:c],[x_1:y_1:z_1],\dots,[x_n:y_n:z_n])
  = \ord_p(x_1\dotsm x_ny_1\dotsm y_nz_1\dotsm z_n)\log p\;,$$
respectively, for places $v$ corresponding to rational primes $p$.
We use these Weil functions to define $N(D,P^{*}_{a,b,c})$
and $N(E,\phi(P^{*}_{a,b,c}))$.

With these choices, we show that (\03.5.1) holds without the $O(1)$ term.
Each term is a sum over all rational primes $p$, so it will
suffice to show that the inequality holds for each $p$.  Fix a prime $p$.
If $p\nmid abc$, then $p$ does not contribute anything to any of these terms,
so the inequality holds for $p$.  Otherwise, by symmetry we may assume that
$p\mid a$; hence $p\nmid b$ and $p\nmid c$.  The contribution to the left-hand
side at $p$ is then
$$\cases \frac{\ord_p a}{n}\log p & \text{if $n\mid\ord_p a$;} \\
  \left(\fracwithdelims[]{\ord_p a}{n} + 1\right)\log p
    & \text{if $n\nmid\ord_p a$.} \endcases$$
The contribution on the right-hand side is
$\left(1+\frac{\ord_p a}{n}\right)\log p$; hence the inequality holds for
the terms at $p$ in this case, too.\qed
\enddemo

\cor{\03.6}  $N(D,P^{*}_{a,b,c}) \le n(abc) + \frac3n h(\phi(P^{*}_{a,b,c}))
+ O(1)$.
\endit

\demo{Proof}  This is immediate from the definitions of the height and counting
functions on $\Bbb P^2$.\qed
\enddemo

At this point we note that the group $G:=\Gm^{2n-2}$ acts faithfully
on the pair $(\Gamma_n,D)$, via
$$\split & (u_2,\dots,u_n,v_2,\dots,v_n)
  \cdot ([a:b:c],[x_1:y_1:z_1],\dots,[x_n:y_n:z_n]) \\
  &\qquad = ([a:b:c],
    [x_1u_2^{-2}u_3^{-3}\dotsm u_n^{-n}:y_1v_2^{-2}\dotsm v_n^{-n}:z_1], \\
  &\qquad\qquad [x_2u_2:y_2v_2:z_2],\dots,[x_nu_n:y_nv_n:z_n])\;.
  \endsplit\tag\03.8$$

One would like to apply Conjecture \ConjB{} to the pair $(\Gamma_n,D)$,
but it is not regular.  However, we do have the following.

\lemma{\03.7}  The pair $(\Gamma_n,D)$ is regular outside of the set
$$x_1\dotsm x_n = y_1\dotsm y_n = z_1\dotsm z_n = 0\;.\tag\03.7.1$$
\endit

\demo{Proof}  Let $P$ be a point on $\Gamma_n$ with multihomogeneous coordinates
$$([a:b:c],[x_1:y_1:z_1],\dots,[x_n:y_n:z_n])\;.$$
It is a regular point of $\Gamma_n$ if and only if the matrix
$$\spreadmatrixlines{1\jot}
  \bmatrix 0 & 0 & 0 &
    \frac{\prod x_i^i}{x_1}\! & \dots & \!\frac{n\prod x_i^i}{x_n} &
    \frac{\prod y_i^i}{y_1}\! & \dots & \!\frac{n\prod y_i^i}{y_n} &
    \frac{\prod z_i^i}{z_1}\! & \dots & \!\frac{n\prod z_i^i}{z_n} \\
  1 & 1 & 1 & 0 & \dots & 0 & 0 & \dots & 0 & 0 & \dots & 0 \\
  \prod y_i^i & \!-\prod x_i^i & 0 &
    \frac{-b\prod x_i^i}{x_1}\! & \dots & \!\frac{-nb\prod x_i^i}{x_n} &
    \frac{a\prod y_i^i}{y_1}\! & \dots & \!\frac{na\prod y_i^i}{y_n} &
    0 & \dots & 0 \endbmatrix$$
has rank $3$.  Assume that (\03.7.1) does not hold at $P$.  By symmetry,
we may assume that $x_1\dotsm x_n\ne0$.  Then the second, third, and fourth
columns are linearly independent, so $\Gamma_n$ is regular at $P$.

It remains to be shown that $D$ has normal crossings at such points $P$.
This only needs to be checked if $P\in\Supp D$.  Again assume that
$x_1\dotsm x_n\ne0$; since $P\in\Supp D$ we may assume by symmetry
that $z_1\dotsm z_n=0$.  Then, by (\03.3a), $y_1\dotsm y_n\ne0$.
Then the rational map $\Gamma_n\dashrightarrow\Bbb A^{2n-1}$ defined by
$$([a:b:c],[x_1:y_1:z_1],\dots,[x_n:y_n:z_n])
  \mapsto \left(\frac{x_2}{y_2},\dots,\frac{x_n}{y_n},
    \frac{z_1}{y_1},\dots,\frac{z_n}{y_n}\right)$$
is regular at $P$.  Then the first, third, and fourth columns of the matrix
mentioned earlier are linearly independent, implying that the above map
is \'etale at $P$.  Thus $D$ has normal crossings at $P$.\qed
\enddemo

The following lemma constructs a well-behaved resolution of the above
singularities.

\lemma{\03.9}  There exist a complete variety $\Gamma_n'$, a divisor $D'$
on $\Gamma_n'$, and a birational morphism $\psi\:\Gamma_n'\to\Gamma_n$
such that:
\roster
\myitem i.  $\Supp D'=\psi^{-1}(\Supp D)$,
\myitem ii.  the pair $(\Gamma_n',D')$ is regular,
\myitem iii.  $\psi$ is an isomorphism over the set
$$\bigl\{P\in\Gamma_n:\vtop{\hbox{$(\Gamma_n,D)$ is regular at $P$ and}
  \hbox{$\phi(P)\notin\Supp E\bigr\}\;,$}}\tag\03.9.1$$
\myitem iv.  $\psi^{*}\phi^{*}\Cal O(1) \le \Cal K_{\Gamma_n'} + D'$
relative to the cone of effective divisors on $\Gamma_n'$, and
\myitem v.  the group action (\03.8) extends to an action on $(\Gamma_n',D')$.
\endroster
\endit

\demo{Proof}  We start by noting that if $[a:b:c]$ and $[a':b':c']$ are
points on $\Bbb P^2$ with $a+b+c=a'+b'+c'=0$ and
$[a:b:c],[a',b',c'] \notin \Supp E\;$, then the morphism
$$\split& ([a:b:c],[x_1:y_1:z_1],\dots,[x_n:y_n:z_n]) \\
  &\qquad \mapsto ([a':b':c'],[a'bcx_1:ab'cy_1:abc'z_1],
    [x_2:y_2:z_2],\dots,[x_n:y_n:z_n])\endsplit$$
induces an isomorphism $\phi^{-1}([a:b:c])\cong\phi^{-1}([a':b':c'])$.
Moreover, this isomorphism preserves the restriction of $D$ to the fibers,
as well as the group action.

Let $F$ be one such fiber, and let $D_F$ denote the restriction of $D$ to $F$.
Also, identify $\Bbb P^1$ with the line $a+b+c=0$ in $\Bbb P^2$.
Then the above isomorphisms define a rational map $\Gamma_n\to F\times\Bbb P^1$.
This rational map is an isomorphism away from the set $\phi^{-1}(\Supp E)$.
Moreover, it is $G$\snug-equivariant.

Let $\rho\:F^{*}\to F$ be an embedded resolution of $D_F$ on $F$.
By \cite{B-M} we may choose $\rho$ such that the action of $G$ on $F$
extends to $F^{*}$ and such that $\rho$ is an isomorphism over the set
where $(F,D_F)$ is regular.  This defines a rational map
$\Gamma_n\dashrightarrow F^{*}\times\Bbb P^1$.  Let $\Gamma_n^{*}$
be the closure of the graph of this map, with projections
$\alpha\:\Gamma_n^{*}\to\Gamma_n$
and $\tau\:\Gamma_n^{*}\to F^{*}\times\Bbb P^1$.
Then $\alpha$ is an isomorphism over the set (\03.9.1), and $\tau$ is
an isomorphism away from $(\phi\circ\alpha)^{-1}(\Supp E)$.

Finally, let $\beta\:\Gamma_n'\to\Gamma_n^{*}$ be an embedded resolution
of $\alpha^{-1}(\Supp D)$ in $\Gamma_n^{*}$, let $D'$ be the resulting normal
crossings divisor on $\Gamma_n'$, and let $\psi=\alpha\circ\beta$.
Again, we may choose $\beta$ so that the group action extends and
such that $\beta$ is an isomorphism over the open set where
$(\Gamma_n^{*},\alpha^{-1}(\Supp D))$ is regular..
Thus the pair $(\Gamma_n',D')$ satisfies conditions (i), (ii), and (v).
Condition (iii) also holds, since the same is true for $\alpha$ and $\beta$.
$$\CD \Gamma_n' \\
  @V\beta VV \\
  \Gamma_n^{*} @>\tau>> F^{*}\times\Bbb P^1 \\
  @V\alpha VV @VV \rho\times\Id_{\Bbb P^1} V \\
  \Gamma_n &\dashrightarrow& F\times\Bbb P^1\endCD$$

Finally, consider Condition (iv).  To begin, note that if $a=0$ at some point
$P\in\Gamma_n$, then (\03.3a) implies that $\prod x_i=0$, so $P\in\Supp D$.
Let $D_{F^{*}}$ be the normal crossings divisor lying over $D_F$.  Then
$$\Supp D'=(\tau\circ\beta)^{-1}
  \bigl((\Supp D_{F^{*}}\times\Bbb P^1) \cup (F^{*}\times\Supp E)\bigr)\;.$$
Since $F^{*}$ is a toric variety with principal orbit
$F^{*}\setminus D_{F^{*}}$, Lemma \02.9 implies that its
logarithmic canonical line sheaf $\Cal K_{F^{*}}(D_{F^{*}})$ is trivial;
hence the logarithmic canonical divisor of
$(F^{*}\times\Bbb P^1)\setminus(D_{F^{*}}\times\Bbb P^1 + F^{*}\times E)$
is the pull-back of $\Cal O(1)$ from the second factor.
Thus the line sheaf $\Cal K_{\Gamma_n'}(D')\otimes\psi^{*}\phi^{*}\Cal O(-1)$
is the line sheaf associated to the logarithmic ramification divisor of
$\tau\circ\beta$.  This divisor is effective, so Condition (iv) holds.\qed
\enddemo

By construction and Lemma \03.7, $\psi$ is an isomorphism over the point
$P_{a,b,c}^{*}\in\Gamma_n$ for all triples $a,b,c$ of relatively prime
integers with $a+b+c=0$, so there is a well-defined point
$P'_{a,b,c}=\psi^{-1}(P^{*}_{a,b,c})$ in $\Gamma_n'$.

\lemma{\03.10}  For relatively prime integers $a,b,c$ with $a+b+c=0$,
$$N(D',P'_{a,b,c}) \le N(D,P^{*}_{a,b,c}) + O(1)\;.$$
\endit

\demo{Proof}  Since $D'$ is reduced and $\Supp D'=\psi^{-1}(\Supp D)$,
the divisor $\psi^{*}D-D'$ is effective.  The inequality then follows from
additivity and functoriality properties of the counting function.\qed
\enddemo

\lemma{\03.11}  For relatively prime integers $a,b,c$ with $a+b+c=0$,
$$h_{\Cal O(1,1,\dots,1)}(P^{*}_{a,b,c}) \le 4h([a:b:c]) + O(1)\;.$$
\endit

\demo{Proof}  By construction
$$P^{*}_{a,b,c} = ([a:b:c],[x_1:y_1:z_1],\dots,[x_n:y_n:z_n])\;,$$
where the $x_i$, $y_i$, and $z_i$ are as in the construction of $P_{a,b,c}$.
In particular, $x_i,y_i,z_i$ are triples of relatively prime integers
for each $i$.  Thus, we may take
$$h_{\Cal O(1,1,\dots,1)}(P^{*}_{a,b,c}) = \log\max\{|a|,|b|,|c|\}
  + \sum_{i=1}^n\log\max\{|x_i|,|y_i|,|z_i|\}\;.$$
The first term equals $h([a:b:c])$, and it is easy to see from
the construction that the second term is bounded from above by
$3h([a:b:c])$.\qed
\enddemo

We are now ready to state and prove the main theorem of the section.

\thm{\03.12}  If, for some $n$, Conjecture \ConjB{} holds for the pair
$(\Gamma_n',D')$ over $\Bbb Q$, then the abc conjecture holds for
all $\epsilon>3/n$.
\endit

\demo{Proof}  Let $\Cal A=\psi^{*}\Cal O(1,1,\dots,1)$.  Since it is the
pull-back of a big line sheaf by a birational morphism, it is big.
Pick $\epsilon'$ such that $4\epsilon'\le \epsilon-3/n$.
Then, by (\03.4), Condition (iv) of Lemma \03.9, properties of
heights and proximity and counting functions, the assumption that Conjecture
\ConjB{} holds, Lemma \03.10, Corollary \03.6, and Lemma \03.11, we have:
$$\split h([a:b:c])
  &= h_{\psi^{*}\phi^{*}\Cal O(1)}(P'_{a,b,c}) + O(1) \\
  &\le h_{\Cal K_{\Gamma_n'}(D')}(P'_{a,b,c}) + O(1) \\
  &= m(D',P'_{a,b,c}) + h_{\Cal K_{\Gamma_n'}}(P'_{a,b,c})
    + N(D',P'_{a,b,c}) + O(1) \\
  &\le \epsilon'\,h_{\Cal A}(P'_{a,b,c})
    + N(D,P^{*}_{a,b,c}) + O(1) \\
  &\le \epsilon'\,h_{\Cal A}(P'_{a,b,c})
    + n(abc) + \frac3n h(P^{*}_{a,b,c}) + O(1) \\
  &\le n(abc) + \epsilon\,h([a:b:c]) + O(1)\endsplit$$
for all $a,b,c$ such that $P'_{a,b,c}$ lies outside of a certain proper
Zariski-closed subset $Z$ of $\Gamma_n'$.  This set is the union of the
exceptional set for Conjecture \ConjB{} and the base locus of the
effective divisor implicit in Lemma \03.9(iv).  As in Section \02,
we may reduce to the case where $Z$ is $G$\snug-invariant, so
(after ignoring components contained in $\Supp D'$) it is a finite union
of fibers of $\phi\circ\psi\:\Gamma_n'\to\Bbb P^1$.  Each fiber can contain
at most two points $P'_{a,b,c}$ (corresponding to $a,b,c$ and $-a,-b,-c$),
so the exceptional set can be eliminated by adjusting the $O(1)$ term.\qed
\enddemo

\cor{\03.13}  If Conjecture \ConjB{} holds for all pairs $(X,D)$ over $\Bbb Q$,
then the abc conjecture holds.
\endit

\beginsection{\04}{Conjecture \ConjB{} in the Split Function Field Case}

This section shows that, in the split function field case of characteristic
zero, the counterpart to Conjecture \ConjB{} holds.  We begin with a
definition.

\defn{\04.1}  Let $C$ be a smooth projective curve over an algebraically
closed field, let $S$ be a subset of $C$, and let
$D=\sum_{P\in C}n_P\cdot[P]$ be a divisor on $C$.
Then
$$\deg_S D = \sum_{P\in S}n_P\;.$$
\endit

\thm{\04.2}  Let $k_0$ be an algebraically closed field, and let $(X,D)$ be a
regular pair with $X$ complete over $k_0$.  Assume that a semiabelian variety
$G$ of dimension $\dim X-1$
acts faithfully on $(X,D)$.  Let $C$ be a smooth projective curve over $k_0$
and let $S$ be a finite set of points on $C$.  Then there exists a proper
Zariski-closed subset $Z$ of $X$ such that any map $f\:C\to X$
with $\Image f\nsubseteq Z\cup\Supp D$ satisfies the inequality
$$\deg_{C\setminus S}(f^{*}D)_{\text{red}}
  \ge \deg f^{*}\Cal K_X(D) - \max\{0,2g(C)-2+\#S\}\tag\04.2.1$$
and therefore the inequality
$$\deg f^{*}\Cal K_X + \deg_S f^{*}D \le \max\{0,2g(C)-2+\#S\}\;.
  \tag\04.2.2$$
\endit

The left-hand side of the first inequality is the counterpart, in the split
function field case, to the {\it truncated counting function\/}
$N^{(1)}(D,\cdot)$.  Conjecture \ConjA{} has also been posed for
truncated counting functions \cite{Vo~2}.

\demo{Proof}  As was noted in Section \02, we may assume that the map
$\pi\:U\to V$ of Lemma \02.3 extends to a morphism $\bar\pi\:X\to\widebar V$
for some projective completion $\widebar V$ of $V$, and that $\Supp D$
contains all degenerate orbits.  In this case these changes do not weaken
the inequality at all, even up to $O(1)$.  Let $Z$ be as in Proposition \02.10.

It will suffice to prove (\04.2.1), since (\04.2.2) follows immediately
by applying the trivial inequality
$$\deg_{C\setminus S}(f^{*}D)_{\text{red}} \le \deg_{C\setminus S}f^{*}D
  = \deg f^{*}D - \deg_S f^{*}D\;.$$

Consider the diagram
$$\CD C @>f>> X @>(0,\Id_X)>> G\times X @>\alpha>> X\times X @>\pr_1>> X \\
  &&&& @VV\pr_2V @VV\pr_2V \\
  &&&& X &=& X \endCD$$
in which $\alpha$ is the map $(\gamma,x)\mapsto(\gamma x,x)$.  This induces
a map of relative tangent bundles:
$$ \CD T_{G\times X/X} @>\alpha_{*}>> \alpha^{*}T_{X\times X/X}
  &\;\cong& \alpha^{*}\pr_1^{*}T_X \\
  &&&& \cup \\
  &&&& \alpha^{*}\pr_1^{*}\bigl(T_X(-\log D)\bigr) \;. \endCD$$
This is a map of sheaves on $G\times X$.
Since the action of $G$ preserves $D$, the image actually lies in
$\alpha^{*}\pr_1^{*}\bigl(T_X(-\log D)\bigr)$.  Let $n=\dim X$.
Taking $\wedge^{n-1}$ of everything gives:
$$\Cal O_{G\times X} \cong \wedge^{n-1}T_{G\times X}/X
  @>\alpha_{*}>> \alpha^{*}\pr_1^{*}\wedge^{n-1}T_X(-\log D)\;.$$
Pulling back via the map $(0,\Id_X)\:X\to G\times X$ gives:
$$\split \Cal O_X
  &@>>> (0,\Id_X)^{*}\alpha^{*}\pr_1^{*}\wedge^{n-1} T_X(-\log D) \\
  &= \wedge^{n-1}T_X(-\log D) \\
  &\cong \Omega^1_X(\log D)\otimes\bigl(\Cal K_X(D)\bigr)\spcheck\;. \endsplit$$
Since $G$ acts faithfully, the above map is nonzero, so it determines a
nonzero global section
$$\omega
  \in \Gamma\bigl(X, \Omega^1_X(\log D)\otimes(\Cal K_X(D))\spcheck\bigr)\;.$$

The above construction is essentially the same as in Section \02,
except that the divisor $D$ has been added.

Now consider $f^{*}\omega$.

\medbreak

{\bc Case I.}  If $f^{*}\omega=0$ then $f$ is everywhere tangent to
an orbit under the $G$\snug-action; hence $\Image f$ is contained
in the closure of some fixed orbit $Y$.  By assumption $Y$ is a non-degenerate
orbit.  By definition of $Z$, $Y$ lies on a fiber of good reduction
for $\bar\pi$, so the closure $\widebar Y$ is the whole fiber,
the pair $(\widebar Y,D\restrictedto{\widebar Y})$ is regular, and the
adjunction formula gives
$\Cal K_{\widebar Y}\cong\Cal K_X\restrictedto{\widebar Y}$.
But $\widebar Y$ is an equivariant completion of a semiabelian variety
$\widebar Y\setminus D\restrictedto{\widebar Y}$, so by Lemma \02.9,
$$\Cal K_X(D)\restrictedto{\widebar Y}
  \cong \Cal K_{\widebar Y}\bigl(D\restrictedto{\widebar Y}\bigr)
  \cong \Cal O_{\widebar Y}\;;$$
hence $f^{*}\Cal K_X(D)$ is trivial.  This gives (\04.2.1) since, on the
left-hand side, $f^{*}D$ is effective.

{\bc Case II.}  Assume $f^{*}\omega\ne0$.  Then it determines a nonzero
global section of
$$\Cal K_C(\log f^{*}D)\otimes f^{*}(\Cal K_X(D))\spcheck\;.$$
Taking degrees, this implies that
$$2g(C)-2+\deg(f^{*}D)_{\text{red}} - \deg f^{*}\Cal K_X(D) \ge 0\;.$$
Applying the inequality
$$\#S \ge \deg_S(f^{*}D)_{\text{red}}$$
then gives (\04.2.1).\qed
\enddemo

\beginsection{\05}{Conjecture \ConjB{} for holomorphic curves}

This section shows that Conjecture \ConjB{} holds in the case of
holomorphic curves (\ie, Nevanlinna theory).  This is done by methods
analogous to those of Section \04.  This relies on a result of McQuillan
to replace the simple argument based on comparing degrees of line sheaves.

\prop{\05.1}  Let $(X,D)$ be a regular pair with $X$ complete over $\Bbb C$,
let $f\:\Bbb C\to X$ be a holomorphic curve not lying entirely in $\Supp D$,
let $f'\:\Bbb C\to\Bbb P(\Omega^1_X(\log D))$ be its derivative, let
$\Cal O(1)$ be the tautological line sheaf on $\Bbb P(\Omega^1_X(\log D))$,
and let $\Cal A$ be a line sheaf on $X$ whose restriction to the Zariski
closure of the image of $f$ is big.  Then
$$T_{\Cal O(1),f'}(r) \le_\exc N^{(1)}_f(D,r)
  + O(\log^{+} T_{\Cal A,f}(r)) + o(\log r)\;.$$
\endit

Here the notation $\le_\exc$ means that the inequality holds outside of
a set of finite Lebesgue measure, and $N^{(1)}_f(D,r)$ denotes the truncated
counting function.

When $D=0$, this was proved by McQuillan \cite{McQ~2, Theorem A}.
A proof of the general case is to appear in \cite{McQ~3}.  Appendix \0A of
this paper gives a more classical proof, based on the lemma on the
logarithmic derivative, using ideas from P.-M. Wong \cite{W}.

\cor{\05.2}  Let $X$, $D$, $f$, $f'$, and $\Cal A$ be as above,
let $\Cal L$ be a line sheaf on $X$, let $d\in\Bbb Z_{>0}$, and let $\omega$
be a global section of $S^d\bigl(\Omega^1_X(\log D)\bigr)\otimes\Cal L\spcheck$.
If $f^{*}\omega\ne0$, then
$$T_{\Cal L,f}(r) \le_\exc d\cdot N^{(1)}_f(D,r)
  + O(\log^{+} T_{\Cal A,f}(r)) + o(\log r)\;.$$
\endit

\demo{Proof}  Let $p\:\Bbb P(\Omega^1_X(\log D))\to X$ be the canonical
projection.  The section $\omega$ corresponds to a global section
$$\omega' \in\Gamma\bigl(\Bbb P(\Omega^1_X(\log D)),
  \Cal O(d)\otimes p^{*}\Cal L\spcheck\bigr)\;,$$
and $f^{\prime*}\omega'=f^{*}\omega$.  Thus the image of $f'$
is not contained in the base locus of $\Cal O(d)\otimes p^{*}\Cal L\spcheck$, so
$$T_{\Cal O(d),f'}(r) \ge T_{\Cal L,f}(r) + O(1)\;.\qed$$
\enddemo

The above corollary now makes it easy to prove the following counterpart to
Conjecture \ConjB{} in Nevanlinna theory.

\thm{\05.3}  Let $(X,D)$ be a regular pair with $X$ complete over $\Bbb C$.
Assume that a semiabelian variety $G$ of dimension $\dim X-1$ acts faithfully
on $(X,D)$.  Let $\Cal A$ be a big line sheaf on $X$.
Then there is a proper Zariski-closed subset $Z\subseteq X$ such that
any holomorphic curve $f\:\Bbb C\to X$ with $\Image f\nsubseteq Z\cup\Supp D$
satisfies the inequality
$$N^{(1)}_f(D,r) \ge_\exc T_{\Cal K_X(D),f}(r)
  - O(\log^{+} T_{\Cal A,f}(r)) - o(\log r)\tag\05.3.1$$
and therefore also the inequality
$$T_{\Cal K_X,f}(r) + m_f(D,r) \le_\exc
  O(\log^{+} T_{\Cal A,f}(r)) + o(\log r)\tag\05.3.2$$
\endit

\demo{Proof}  As in Section \04, we may assume that the map
$\pi\:U\to V$ of Lemma \02.3 extends to a morphism $\bar\pi\:X\to\widebar V$
for some projective completion $\widebar V$ of $V$, and that $\Supp D$
contains all degenerate orbits.  These changes weaken the inequalities only
up to $O(1)$.  Let $Z_0$ be as in $Z$ of Proposition \02.10.

Let $\pi\:X'\to X$ be a birational morphism, with $X'$ projective.
By Kodaira's lemma (see, for example, \cite{Vo~1, Prop.~1.2.7}) some
positive tensor power of $\pi^{*}\Cal A$ is isomorphic to
$\Cal A'\otimes\Cal O(D)$, with $\Cal A'$ ample and $D$ effective.
Let $Z_1$ be the image under $\pi$ of the base locus of $D$.
Then let $Z=Z_0\cup Z_1$, and assume that the image of $f$ is not
contained in $Z\cup\Supp D$.

Also as in Section \04, it will suffice to prove (\05.3.1).

Let $\omega$ be the form constructed in the proof of Theorem \04.2.
As in that proof, if $f^{*}\omega=0$, then the image of $\omega$ is
contained in the closure of a non-degenerate orbit $Y$ of the group action,
and (\05.3.1) holds again since $f^{*}\Cal K_X(D)$ is trivial.

Otherwise, $f^{*}\omega\ne0$, and (\05.3.1) follows immediately from
Corollary \05.2.\qed
\enddemo

\beginsection{\06}{A hierarchy of problem types}

For several decades, it has been known that valuable insight into a
diophantine problem over number fields can be gained by looking at the
corresponding problem over function fields.  In the function field case,
one looks first at the split case (where everything is defined over the
field of constants of the function field).

More recently, it has been observed that diophantine problems are formally
similar to problems in Nevanlinna theory, and that insight into the former
may be gained by looking at the latter (and sometimes vice versa).

This section introduces a hierarchy of problem types that incorporates
the above observations, plus a few others.

Given a classification of regular pairs $(X,D)$, one can pose a number
of related problems in various contexts:
\roster
\item"$\bullet$" Find the (algebraic) {\bc exceptional set}; i.e., the Zariski
closure of the union of the images of all non-constant strictly rational maps
$G\dashrightarrow X\setminus D$, where $G$ is either $\Gm$
or an abelian variety.  A {\bc strictly rational map} \cite{I, \S\kern.2ex2.12}
is a rational map $X\dashrightarrow Y$ such that the closure of the graph
is proper over $X$.
\item"$\bullet$" For each $\epsilon>0$, find the exceptional subset $Z$
for the main inequality of Conjecture \ConjA.
\item"$\bullet$" Prove the inequality of Conjecture \ConjA{} in the split
function field case of characteristic zero.
\item"$\bullet$" Prove that all non-constant holomorphic curves
$\Bbb C\to X\setminus D$ must lie in the exceptional set.
\item"$\bullet$" Prove the inequality of Conjecture \ConjA{}
for holomorphic curves $\Bbb C\to X$.
\item"$\bullet$" Prove that the set of integral points on $X\setminus D$ is not
Zariski-dense, in the (general) function field case of characteristic zero.
\item"$\bullet$" Prove the inequality of Conjecture \ConjA{}
in the function field case of characteristic zero.
\item"$\bullet$" Prove that the set of integral points on $X\setminus D$ is not
Zariski-dense, in the function field case of characteristic $p>0$.
\item"$\bullet$" Prove the inequality of Conjecture \ConjA{}
in the function field case of characteristic $p>0$.
\item"$\bullet$" Prove that the set of integral points on $X\setminus D$ is not
Zariski-dense, in the number field case.
\item"$\bullet$" Prove the inequality of Conjecture \ConjA{}
in the number field case.
\item"$\bullet$" Prove ``moving targets'' versions of the above.
\endroster

For example, one may pose Conjecture \ConjA{} in each of the above contexts
(in which case the first, fourth, sixth, eighth, ninth, and tenth entries
would not apply, since they are false in some cases).  If one restricted
to pairs $(X,D)$ of
logarithmic general type (\ie, $\Cal K_X(D)$ is big), then all of the above
would apply.  Or, one may pose other restrictions, as was done for
Conjecture \ConjB.

The above hierarchy is ranked roughly from easiest to hardest.  The general
idea is that one would start from the top and work down from there, using
the insight gained on earlier steps to help with the later ones.  For example,
this paper works through some of the above steps for Conjecture \ConjB{}.
The first item is not useful for Conjecture \ConjB{} (since the exceptional
set is all of $X$), the second item is solved in Section \02, the third
in Section \04, the fourth item again is not useful, and the fifth item was
solved in Section \05.  The sixth item is not useful, so the next step is
to try to prove Conjecture \ConjB{} in the function field case of
characteristic zero.

As another example, McQuillan's paper \cite{McQ~2} proved the Nevanlinna-theory
analogue of Bogomolov's theorem bounding the number of curves of given
genus on surfaces of general type with $c_1^2>c_2$.  This may be regarded
as proceeding from the second step to the fourth.

\beginsection{\07}{Complements}

This section gives some variations of the method of Section \03.
These give some interesting implications.

We begin with a variation that shows that Conjecture \ConjB, with $D=0$,
would still imply a weak form of the abc conjecture, with (\00.1.1)
replaced by
$$\max\{|a|,|b|,|c|\} \le C\prod_{p\mid abc}p^{24+\epsilon}\;.\tag\07.1$$
This shows that Conjecture \ConjB{} leads to nontrivial results even on
rational varieties without a divisor.  This does not augur well for the
prospects of actually proving Conjecture \ConjB{} by known methods, since
the absence of $D$ and the triviality of the Albanese rule out the
usual gains to be expected from taking a carefully chosen line sheaf
on a product of several copies of the variety.

\prop{\07.2}  There is a rational three-fold $X$ with the property that
if Conjecture \ConjB{} holds for $X$ (with $D=0$) over $\Bbb Q$, then
Conjecture \00.1 holds with (\00.1.1) replaced by (\07.1).
\endit

\demo{Proof}  For relatively prime integers $a,b,c$ with $a+b+c=0$,
pick integers $u$, $v$, $w$, $x$, $y$, and $z$ such that $x$, $y$, and $z$
are as large as possible, and such that
$$ux^5=a,\qquad vy^5=b,\qquad\text{and}\qquad wz^5=c\;.\tag\07.2.1$$
This defines a point
$$P_{a,b,c} = \bigl([u:v:w],[x:y:z]\bigr)$$
on the subvariety $X$ of $\Bbb P^2\times\Bbb P^2$ cut out by the equation
$$ux^5+vy^5+wz^5=0\;.$$
This variety is smooth and rational, since the projection to the second factor
exhibits it as a $\Bbb P^1$\snug-bundle over $\Bbb P^2$.  It admits an action
of $\Gm^2$ by
$$(g_1,g_2) \cdot ([u:v:w],[x:y:z])
  = ([ug_1^{-5}:vg_2^{-5}:w],[xg_1:yg_2:z])\;.$$

For $i=1,2$ let $\pr_i\:X\to\Bbb P^2$ denote the projection to the
$i^{\text{th}}$ factor.  Then, by the adjunction formula,
$\Cal K_X\cong\pr_1^{*}\Cal O(-2)\otimes\pr_2^{*}\Cal O(2)$.
By functoriality of heights and by Conjecture \ConjB, it follows that
$$2h([x:y:z]) - 2h([u:v:w]) \le \epsilon(h([x:y:z])+h([u:v:w])) + O(1)$$
outside a proper Zariski-closed subset $Z\subsetneq X$.  Therefore,
$$h([x:y:z]) \le h([u:v:w]) + \epsilon\,h([a:b:c]) + O(1)\tag\07.2.2$$
(with a different $\epsilon$).  Thus, by (\07.2.1) and (\07.2.2),
$$\split h([a:b:c]) &\le h([u:v:w]) + 5h([x:y:z]) \\
  &\le 6h([u:v:w]) + 5\epsilon h([a:b:c]) + O(1) \\
  &\le 24n(abc) + 5\epsilon h([a:b:c]) + O(1)\;.\endsplit$$
Here the last step follows because any given prime may occur in $uvw$
to at most a fourth power.  This gives (\07.1) (with a different $\epsilon$).

As in Section \03, the set $Z$ corresponds to at most finitely many triples
$a,b,c$, and so it can be ignored.\qed
\enddemo

\remk{\07.3}  The above approach will also work if the exponents are changed
from $5$ to $4$; this gives $27+\epsilon$ instead of $24+\epsilon$.
\endit

\remk{\07.4}  The above variety $X$ (modified as in Remark \07.3) also appears
in conjunction with the abc conjecture in the paper \cite{McQ~1},
where it is noted that a different conjecture also would imply the
abc conjecture.
\endit

We mention here one more example of this sort, involving in this case
a surface that can be explicitly described.  We consider Pell's equation
$$x^2-dy^2=\pm4,\qquad d\in\Bbb Z\;.$$
This equation potentially infinitely many solutions $(x,y)\in\Bbb Z^2$.

We consider here the question of whether these solutions are {\bc mostly
square free}; \ie, whether for all fixed $\epsilon>0$, the largest square
factor of $x$ or $y$ is $O(\max\{|x|,|y|\})$.  If the abc conjecture
is true, then an easy argument implies that such solutions are mostly
square free.  Although the converse does not seem to hold, the question
of whether solutions to Pell's equation are mostly square free still captures
some of the flavor of the abc conjecture.

\prop{\07.5}  There is a rational projective surface $X$ and a divisor $D$
on $X$ such that $(X,D)$ is regular, and such that if Conjecture \ConjB{}
holds for $(X,D)$ over all quadratic fields, then solutions of Pell's equation
are mostly square free.
\endit

\demo{Proof}  This follows by applying essentially the same methods as in
Section \03 and the earlier part of this section, applied to the equation
$$x^2v^4-y^2w^4=\pm4\;,$$
which determines a variety in $\Bbb P^1\times\Bbb P^2$ whose desingularization
is easy to find explicitly.  The resulting variety $X$ can be described
pictorially below.
$$
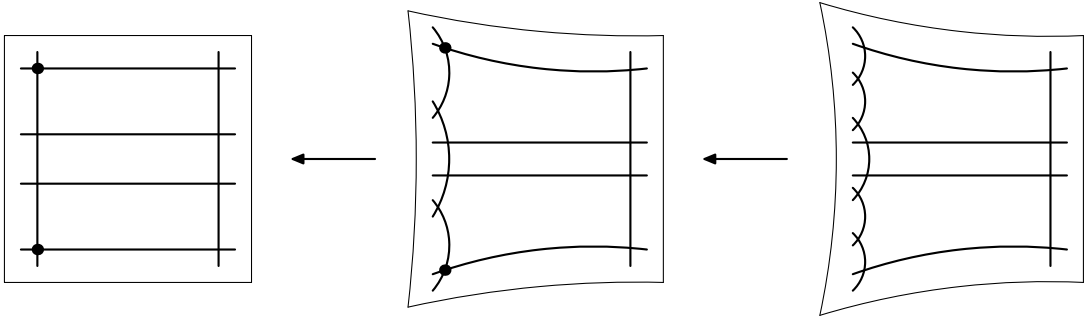
\centerline{\epsfbox{ratabc.mps}}$$
Here the diagram on the left depicts the divisor
$$\pr_1^{*}([0]+[\infty]+[1]+[-1])+\pr_2^{*}([0]+[\infty])$$
in $\Bbb P^1\times\Bbb P^1$, with the lines $\pr_1^{-1}(\pm1)$ drawn in
the middle.  The arrows are blowings-up at the fat points.
Since $\Bbb P^1\times\Bbb P^1$ blown up at one point is isomorphic to the
blowing-up of $\Bbb P^2$ at two points, $X$ can also be described as a certain
blowing-up of $\Bbb P^2$, with $D$ being the inverse image of the coordinate
axes and the lines $y=\pm x$.

We leave the details of this proof to the reader.\qed
\enddemo

Although the surface $X$ is just a few blowings-up away from the more
general version of Schmidt's Subspace Theorem \cite{Vo~1, Thm.~2.5.8},
it is still, unfortunately, out of reach.

\genbeginsect{Appendix \0A.  Proof of Proposition \05.1}

This appendix provides a proof of Proposition \05.1, since a proof has not yet
appeared in print.  The proof presented here relies on a geometric
formulation of the classical ``lemma on the logarithmic derivative'';
see Theorem \0A.2.  This may be of independent interest.

This geometric logarithmic derivative lemma was first proved by
M. McQuillan \cite{McQ~2}, but was stated only implicitly.  See also
\cite{McQ~3}.  Theorem \0A.6 is also due to him; the proof given here
is based on his ideas but uses more elementary methods.
The proof of Theorem \0A.2 presented here does not follow McQuillan;
rather, it reduces to the classical logarithmic derivative lemma by
a method originally due to P.-M. Wong \cite{W}.  In \cite{W, Thm.~4.1},
Wong proved the following special case:  If $X$ is a smooth projective variety,
if $D$ is an effective divisor with strict normal crossings (Definition
\0A.2.3), if $f\:\Bbb C\to X\setminus D$ is a non-constant holomorphic map,
if $\omega$ is a global section of $\Cal J_k^mX(\log D)$, and
if $(j^kf)^{*}\omega$ is not identically zero, then
$$\int_0^{2\pi}\log^{+}\bigl|\omega(j^k f(re^{i\theta}))\bigr|
    \,\frac{d\theta}{2\pi}
  \le_\exc O(\log T_{\Cal A,f}(r)) + O(\log r)\;.$$

We begin with the geometric logarithmic derivative lemma.
Let $X$ be a smooth compact complex algebraic
variety, let $D$ be a normal crossings divisor on $X$, and let $f\:X\to\Bbb C$
be a holomorphic curve whose image is not contained in the support of $D$.
Let $V=\Bbb V(\Omega^1_X(\log D))$ and
$\widebar V=\Bbb P(\Omega^1_X(\log D)\oplus\Cal O_X)$.  Here we use the
conventions of \cite{EGA~II, 1.7.8}, so that $V$ is the total space of
the tangent sheaf of $X$ with logarithmic zeroes along $D$ and $\widebar V$
is the obvious projective completion.  The natural map
$f^{*}\:\Omega^1_X(\log D)\to\Omega^1_{\Bbb C}(\log f^{*}D)$, together with
the map $\Cal O_X\to\Omega^1_{\Bbb C}$ given by $1\mapsto dz$, gives a
holomorphic map $Df\:\Bbb C\to\Bbb P(\Omega^1_X(\log D)\oplus\Cal O_X)$.
(Note that this differs from the map $f'\:\Bbb C\to\Bbb P(\Omega^1_X(\log D))$
defined in Section \05.)  Let $[\infty]$ be the complement of $V$ in
$\widebar V$, and choose a Weil function (or Green function) $g_{[\infty]}$
for this divisor.  We will use the normalization of (\0A.1), below.

One possible way to choose this Weil function is the following.  Choose a
Hermitian metric on $\Omega^1_X(\log D)$; its dual metric on the logarithmic
tangent bundle $T_X(-\log D)$ induces a Weil function for $[\infty]$
by the formula
$$g_{[\infty]}(P) = \log^{+}\|\xi\|\;,\tag\0A.1$$
where $P\in V$ corresponds to an element $\xi$ in the fiber of $T_X(-\log D)$
over the corresponding point of $X$, and $\log^{+}x$ is defined as
$\max\{\log x,0\}$.

The geometric logarithmic derivative lemma can then be stated as follows.

\thm{\0A.2}  Let $(X,D)$ be a regular pair with $X$ complete over $\Bbb C$,
let $f\:\Bbb C\to X$ be a holomorphic curve
whose image is not contained in the support of $D$, and let $\Cal A$ be a line
sheaf on $X$ whose restriction to the Zariski closure of the image of $f$
is big.  Let $Df$ and $m_{Df}([\infty],r)$ be as above.  Then
$$m_{Df}([\infty],r) \le_\exc O(\log^{+} T_{\Cal A,f}(r)) + o(\log r)\;.
  \tag\0A.2.1$$
\endit

\demo{Proof}  The proof works by reducing to a situation where the classical
lemma on the logarithmic derivative can be applied a finite number of times.
We begin with some reductions.

\lemma{\0A.2.2}  Let $(X_2,D_2)$ be a regular pair with $X_2$
complete over $\Bbb C$, and let $\pi\:X_2\to X$ be a morphism that
induces a birational morphism $X_2\to\pi(X_2)$.  Assume also that
$\Supp D_2=\pi^{-1}(\Supp D)$ and that the image of $f$ is a
Zariski-dense subset of $\pi(X_2)$.
Let $g\:\Bbb C\to X_2$ be the (unique) holomorphic curve such that
$\pi\circ g=f$.  Let $\Cal A_2$ be a big line sheaf on $X_2$.
Finally, assume that Theorem \0A.2 holds for $g$ and $\Cal A_2$.
Then it also holds for $f$ and $\Cal A$.
\endit

\demo{Proof}  Let $V_2=\Bbb V(\Omega^1_{X_2}(\log D_2))$ and let
$\widebar V_2=\Bbb P(\Omega^1_{X_2}(\log D_2)\oplus\Cal O_{X_2})$.
The natural map $\pi^{*}\Omega^1_X\to\Omega^1_{X_2}$ induced by $\pi$
extends to a map $\pi^{*}\Omega^1_X(\log D)\to\Omega^1_{X_2}(\log D_2)$.
This defines a morphism $V_2\to V\times_X X_2$ and hence
a morphism $D\pi\:V_2\to V$.
Moreover, $Df\restrictedto{(Dg)^{-1}(\Supp D_2)}
  = D\pi\circ Dg\restrictedto{(Dg)^{-1}(\Supp D_2)}$.  We claim that
$$m_{Df}([\infty],r) \le m_{Dg}([\infty],r) + O(1)\;.\tag\0A.2.2.1$$
Indeed, we may assume that the Weil functions for the respective divisors
$[\infty]$ on $\widebar V$ and $\widebar V_2$ are of the form (\0A.1).
The inequality then follows from the fact that the natural map
$\theta\:T_{X_2}(-\log D_2)\to\pi^{*}T_X(-\log D)$ satisfies
$\|\theta(\xi)\|\ll\|\xi\|$ for all $\xi\in T_{X_2}(-\log D_2)$
by compactness of $X_2$.

To compare the error terms, we have
$$T_{\Cal A,f}(r) \gg T_{\Cal A_2,g}(r) + O(1)\tag\0A.2.2.2$$
since $\pi^{*}\Cal A$ is big.
The lemma then follows from (\0A.2.2.1) and (\0A.2.2.2).\qed
\enddemo

Let $Z$ be the Zariski closure of the image of $f$.  By letting
$\pi\:X_2\to X$ be a resolution of the pair $\bigl(Z,D\restrictedto Z\bigr)$,
we may assume that the image of $f$ is Zariski dense.

By applying Chow's lemma and resolving singularities again, we may further
assume that $X$ is projective.

\defn{\0A.2.3}  A divisor $D$ on a smooth variety $X$ has {\bc strict normal
crossings} (also called {\bc simple normal crossings}) if it is a normal
crossings divisor and all of its irreducible components are regular.
\endit

It is well known that in the present situation there exists a smooth projective
variety $X_2$ and a birational morphism $\pi\:X_2\to X$ such that
$\pi^{-1}(\Supp D)$ is the support of a divisor $D_2$ on $X_2$ with
strict normal crossings, and such that $\pi$ is an isomorphism outside
of the support of $D_2$.  After applying Lemma \0A.2.2 to $\pi$, we may
therefore assume that $D$ has strict normal crossings.

The remainder of the proof follows the method of P.-M. Wong \cite{W}:

Let $n=\dim X$, let $\Cal L_0$ be a very ample line sheaf on $X$, let
$E_{00},\dots,E_{0,2n}$ be effective divisors corresponding to $\Cal L_0$
such that any $n+1$ have empty intersection, and for integers $i,j$ with
$0\le i<j\le 2n$ choose a rational function $f_{ij}$ on $X$ such that
$(f_{ij})=E_{0i}-E_{0j}$.  Then the set
$$\biggl\{\frac{df_{ij}}{f_{ij}}:0\le i<j\le 2n\biggr\}$$
of rational sections of $\Omega^1_X$ has the property that, for each point
$P\in X$, some subset of this set is regular at $P$ and generates $\Omega^1_X$
there.

Next, recalling that $D$ has strict normal crossings, write
$D=D_1+\dots+D_\ell$, where each $D_i$ is effective and has smooth support.
For each $i=1,\dots,\ell$ let $\Cal L_i$ and $\Cal L_i'$ be very ample
line sheaves on $X$ such that $\Cal L_i\cong\Cal L_i'(D_i)$.  For each $i$
and each $j=1,\dots,n$ choose effective divisors $E_{ij}$ and $E_{ij}'$
associated to $\Cal L_i$ and $\Cal L_i'$, respectively, such that
$\bigcap_j\bigl(\Supp E_{ij}\cup\Supp E_{ij}'\bigr)$ is disjoint from
$\Supp D_i$.  For each $i$ and $j$ choose a rational function $g_{ij}$ on $X$
such that $(g_{ij})=D_i+E_{ij}'-E_{ij}$.  Then the set
$$\biggl\{\frac{dg_{ij}}{g_{ij}}:1\le i\le\ell,\;1\le j\le n\biggr\}$$
of rational sections of $\Omega^1_X(\log D)$ has the property that, for
each point $P\in X$, some subset of this set is regular at $P$ and
generates $\Omega^1_X(\log D)$ over $\Omega^1_X$ there.

Let
$$\Cal H = \{f_{ij}:0\le i<j\le 2n\}
  \cup \{g_{ij}:1\le i\le\ell,\;1\le j\le n\}\;.$$
This set has the property that, for each $P\in X$, there is a subset $\Cal H_P$
of $\Cal H$ such that $dh/h$ is a regular section of $\Omega^1_X(\log D)$
at $P$ for all $h\in\Cal H_P$, and these sections generate $\Omega^1_X(\log D)$
there.  By compactness of $X$, it then follows that
$$g_{[\infty]}(Df(z))
  \le \log^{+}\max_{h\in\Cal H} \fracwithdelims||{(h\circ f)'}{h\circ f}
    + O(1) \;.$$
Thus, by the classical lemma on the logarithmic derivative applied to the
meromorphic functions $h\circ f$, $h\in\Cal H$, we have
$$\qed\split m_{Df}([\infty],r)
  & \le \max_{h\in\Cal H} \int_0^{2\pi}\log^{+}
    \fracwithdelims||{(h\circ f)'(re^{i\theta})}{(h\circ f)(re^{i\theta})}
    \frac{d\theta}{2\pi} + O(1) \\
  &\le_\exc \sum_{h\in\Cal H} O(\log^{+} T_{h\circ f}(r)) + o(\log r) \\
  &\le_\exc O(\log^{+} T_{\Cal A,f}(r)) + o(\log r)\;.\endsplit$$
\enddemo

\remk{\0A.3}  When $X=\Bbb P^1$ and $D=[0]+[\infty]$, this theorem reduces
to the classical lemma on the logarithmic derivative.  Indeed, in that case
$T_X(-\log D)$ is the trivial line bundle on $\Bbb P^1$, via the isomorphism
$z\frac\partial{\partial z}\mapsto 1$.  Let $f\:\Bbb C\to\Bbb P^1$ be a
holomorphic map whose image is not contained in $\{0,\infty\}$.  Then
$g_{[\infty]}(Df(z))=\log^{+}|f'(z)/f(z)| + O(1)$, so in this case
$$m_{Df}([\infty],r)
  = \int_0^{2\pi}\log^{+}\fracwithdelims||{f'(re^{i\theta})}{f(re^{i\theta})}
    \frac{d\theta}{2\pi} + O(1)\;,$$
which is the quantity appearing in the classical lemma on the logarithmic
derivative.
\endit

\remk{\0A.4}  Since the above proof merely reduces to multiple applications
of the classical case, any sharpening of the error term in the classical
lemma on the logarithmic derivative leads immediately to a correspondingly
sharp error term in the above generalization (up to a constant factor).
Here we used the error term of \cite{S, Thm.~3.11} or \cite{Y, Thm.~1}.
\endit

Proposition \05.1 is a fairly easy consequence of the above geometric
lemma on the logarithmic derivative.  It will be proved in a slightly stronger
form, involving a modified counting function for ramification.  Recall
that $\widebar V=\Bbb P(\Omega^1_X(\log D)\oplus\Cal O_X)$, and let
$p\:P\to\widebar V$ be the blowing-up of $\widebar V$ along the zero
section (of $V$); \ie, the section corresponding to the projection
$\Omega^1_X(\log D)\oplus\Cal O_X\to\Cal O_X$.  Let $[0]\subseteq P$ be
the exceptional divisor and let $\phi\:\Bbb C\to P$ be the lifting of
$Df\:\Bbb C\to\widebar V$.

\defn{\0A.5}  The {\bc $D$\snug-modified ramification counting function}
of a holomorphic curve $f\:\Bbb C\to X$ is the counting function for
$\phi^{*}[0]$:
$$N_{\text{Ram}(D),f}(r) = N_\phi([0],r)\;.$$
\endit

\thm{\0A.6}  Let $(X,D)$ be a regular pair with $X$ complete over $\Bbb C$,
let $f\:\Bbb C\to X$ be a holomorphic curve not lying entirely within
the support of $D$, and let $\Cal A$ be a line sheaf on $X$ whose restriction
to the Zariski closure of the image of $f$ is big.
Let $f'\:\Bbb C\to \Bbb P(\Omega^1_X(\log D))$ be the canonical lifting
of $f$, and let $\Cal O(1)$ denote the tautological line sheaf on
$\Bbb P(\Omega^1_X(\log D))$.  Then
$$T_{\Cal O(1),f'}(r) \le_\exc N_f^{(1)}(D,r) - N_{\text{Ram}(D),f}(r)
  + O(\log^{+}T_{\Cal A,f}(r)) + o(\log r)\;.\tag\0A.6.1$$
\endit

\demo{Proof}  This proof essentially follows McQuillan \cite{McQ~2}, but
some details are different.

Recall the blowing-up $p\:P\to\Bbb P(\Omega^1_X(\log D)\oplus\Cal O_X)$.
This $P$ admits a morphism $q\:P\to\Bbb P(\Omega^1_X(\log D))$, extending
the rational map
$$\Bbb P(\Omega^1_X(\log D)\oplus\Cal O_X)
  \dashrightarrow \Bbb P(\Omega^1_X(\log D))$$
associated to the canonical map
$\Omega^1_X(\log D)\hookrightarrow \Omega^1_X(\log D)\oplus\Cal O_X$.
We first compare the pullbacks $q^{*}\Cal O(1)$ and $p^{*}\Cal O(1)$
of the tautological line sheaves on $\Bbb P(\Omega^1_X(\log D))$
and $\widebar V=\Bbb P(\Omega^1_X(\log D)\oplus\Cal O_X)$, respectively.
Let $s$ be any nonzero rational section of $\Omega^1_X(\log D)$.  This
determines a rational section of $\Cal O(1)$ on $\Bbb P(\Omega^1_X(\log D))$.
The corresponding divisor $D_1$ is the sum of a component which is generically
a hyperplane section on fibers over $X$, and the pull-back of a divisor on $X$.
But also $(s,0)$ is a nonzero rational section of
$\Omega^1_X(\log D)\oplus\Cal O_X$, giving a rational section of $\Cal O(1)$
on $\widebar V$, hence a divisor $D_2$ which is again the sum of a
generic hyperplane section and the pull-back of a divisor on $X$.
Comparing $q^{*}D_1$ with $p^{*}D_2$, we see that they coincide except that
$p^{*}D_2$ contains $[0]$ with multiplicity $1$.  Hence
$$q^{*}\Cal O(1) \cong p^{*}\Cal O(1) \otimes\Cal O(-[0])\;.$$

The global section $(0,1)$ of $\Omega^1_X(\log D)\oplus\Cal O_X$ corresponds
to the divisor $[\infty]$ on $\widebar V$; hence
$\Cal O([\infty])\cong\Cal O(1)$ and therefore
$$T_{\Cal O(1),f'}(r)
  = N_{Df}([\infty],r) + m_{Df}([\infty],r) - N_\phi([0],r) - m_\phi([0],r)
    + O(1) \;.
  \tag\0A.6.2$$
The lifted curve $Df$ meets $[\infty]$ only over $D$, and with multiplicity
at most $1$.  Hence
$$N_{Df}([\infty],r) \le N^{(1)}_f(D,r)\;.\tag\0A.6.3$$
The second term $m_{Df}([\infty],r)$ is bounded by Theorem \0A.2.
The third term (by definition) is equal to the modified ramified counting
function, and the fourth term on the right-hand side of (\0A.6.2) is bounded
from below and therefore can be ignored.  Thus (\0A.6.2) gives (\0A.6.1).\qed
\enddemo

Finally, we note that Theorem \0A.6 implies Proposition \05.1 since
$N_{\text{Ram}(D),f}(r)\ge0$.

Again, sharper error terms in the classical lemma on the logarithmic derivative
lead to sharper error terms in Theorem \0A.6.

Theorem \0A.6 can also be proved in the context of coverings
(classically called ``algebroid functions'').  In this context, one has
a finite ramified covering $p\:Y\to\Bbb C$ and a holomorphic function
$f\:Y\to X$.  In place of the set $\{z\in\Bbb C:|z|=r\}$,
let $Y\langle r\rangle=\{z\in Y : |p(z)|=r\}$.  The classical lemma on the
logarithmic derivative, originally due to Valiron \cite{Va}, can then
be stated as follows.

\thm{\0A.7}  Let $Y$, $p$, and $f$ be as above.  Then, for all $\epsilon>0$,
$$\int_{Y\langle r\rangle}
    \log^{+}\left|\frac 1{f(z)} \frac{df}{dp}(z)\right| d^c\log|p(z)|^2
  \le_\exc (\deg p+\epsilon)[7\log T_f(r) + 6\log^{+}r]\;.$$
\endit

\demo{Proof}  This follows from \cite{A, Thm.~2.2} with $\tau=|p|^2$,
$\Bbb B=1$, $\Theta=dp$, and other notation as in \cite{A, \S\,1.1}.\qed
\enddemo

To state Theorem \0A.2 in the context of coverings, we define $Df$ to be
the map $Y\to\widebar V$
associated to the map $f^{*}\:\Omega^1_X(\log D)\to\Omega^1_Y(\log f^{*}D)$,
together with the map $\Cal O_X\to\Omega^1_Y$ given by $1\mapsto dp$.
The statement and proof of Theorem \0A.2 then carry over directly,
with (\0A.2.1) replaced by
$$\int_{Y\langle r\rangle} g_{[\infty]}(Df(z)) d^c\log|p(z)|^2
  \le_\exc O(\log T_{\Cal A,f}(r)) + O(\log^{+}r)\;.\tag\0A.8$$
Then Theorem \0A.6 holds with (\0A.6.1) replaced by
$$\split T_{\Cal O(1),f'}(r)
  &\le_\exc N_f^{(1)}(D,r) + N_{\text{Ram},p}(r) - N_{\text{Ram}(D),f}(r) \\
  &\qquad+ O(\log T_{\Cal A,f}(r)) + O(\log^{+}r)\;,\endsplit\tag\0A.9$$
where the additional term $N_{\text{Ram},p}(r)$ is the counting function
for the ramification of $p\:Y\to\Bbb C$.  The proof is again essentially
the same, except that the additional term $N_{\text{Ram},p}(r)$ needs to
be added to the right-hand side of (\0A.6.3).

As has been noted by R. Kobayashi, the above results also hold for higher
jet bundles.  This can be seen by noting that the higher jet bundle is the
tangent bundle of the next lower jet bundle.  The details are left to the
reader.

\enddocument